\documentclass[12pt]{article}

\usepackage{diagbox}
\usepackage{mathtools}
\usepackage{bbm}
\usepackage{latexsym}
\usepackage{epsfig}
\usepackage{amsmath,amsthm,amssymb,enumerate}

\usepackage[a-1b]{pdfx}
\usepackage{hyperref}
\usepackage[normalem]{ulem}

\usepackage{color}

\def\a{\alpha} \def\b{\beta} \def\d{\delta} 
\def\e{\varepsilon} \def\f{\phi} \def\F{{\Phi}}  
  \def\k{\kappa}
\def\z{\zeta} \def\th{\theta}    
 \def\m{\mu}  \def\p{\pi}
   
\def\t{\tau}

\newtheorem{theorem}{Theorem}
\newtheorem{lemma}[theorem]{Lemma}
\newtheorem{corollary}[theorem]{Corollary}

\newcommand{\wh}[1]{\widehat{#1}}

\newcommand{\rdown}[1]{{\left\lfloor #1\right \rfloor}}

\newcommand{\brac}[1]{\left(#1\right)}

\newcommand{\bfrac}[2]{\left(\frac{#1}{#2}\right)}

\def\cE{{\cal E}}

\newcommand{\set}[1]{\left\{#1\right\}}

\def\vK{\vec K}
\def\E{\mathbb{E}}

\def\Pr{\mathbb{P}}

\def\cF{{\cal F}}
\newcommand{\ignore}[1]{}

\def\cE{{\mathcal E}}
\def\cF{{\mathcal F}}

\newcommand{\beq}[2]{\begin{equation}\label{#1}#2\end{equation}}
\newcommand{\mults}[1]{\begin{multline*}#1\end{multline*}}

\usepackage{tikz}
\usetikzlibrary{decorations.pathmorphing}
\usetikzlibrary{positioning}
\usetikzlibrary{arrows,automata}
\usetikzlibrary{shapes.misc}
\usetikzlibrary{backgrounds}
\usetikzlibrary{arrows,shapes}

\def\N{\mathbb{N}}
\def\wC{{\wh C}}

\begin{document}
\author{Patrick Bennett\thanks{Research supported in part by Simons Foundation Grant \#426894.} \and Alan Frieze\thanks{Research supported in part by NSF grant DMS1952285}}

\date{}
\title{Weighted tree games}
\maketitle
\begin{abstract}
We consider a variation on Maker-Breaker games on graphs or digraphs where the edges have random costs. We assume that Maker wishes to choose the edges of a spanning tree, but wishes to minimise his cost. Meanwhile Breaker wants to make Maker's cost as large as possible.
\end{abstract}
\section{Introduction}
We consider a variation on Maker-Breaker games on graphs or digraphs where the edges have costs. We assume that each edge has a (random but known to both players) cost and Maker's goal is build some some structure, e.g. a spanning tree, but wishes to minimise his cost. In a play of the game, Maker chooses an edge for his structure and then Breaker deletes $b$ edges. If Maker is unable to build the desired structure at all we say the cost is infinity. Maker will have to avoid this situation to get a meaningful result. Biased Maker-Breaker games have a long history of research, beginning with Chv\'atal and Erd\H{o}s \cite{CE}. We refer the reader to an excellent monograph of Hefetz, Krivelevich, Stojakovi\'c and Szabo \cite{Book} on this topic.

We begin with a simple case that does not involve graphs:  let $\N=\set{1,2,\ldots,}$ be the set of positive integers. Maker in her turn has to choose an $i\in [n]$ and irrevocably assign a value to $f(i)\in\N$. We let $M_t$ denote the set of elements $i\in[n]$ that have been selected in this way by Maker after $t$ rounds of play. Breaker in his turn selects $i_1,i_2,\ldots,i_b\in [n]$ and $j_1,j_2,\ldots,j_b\in\N$ and makes $j_r$ unavailable to Maker for the value of $f(i_r)$, $r=1,2,\ldots,b$. We let $B_t$ denote the set of pairs $(i,j)$ for which Breaker has made $f(i)=j$ unavailable to Maker after $t$ rounds. Thus $M_0=B_0=\emptyset$ and in round $t$, Maker adds one element to $M_t$ to create $M_{t+1}$ and Breaker adds $b$ pairs to $B_t$ to create $B_{t+1}$.

Let $\f(f)=\sum_{i=1}^nf(i)$. Maker's aim is to keep $\f$ as small as possible and Breaker has the opposite intention. Our first result is the following.
\begin{theorem}\label{th1}
If Breaker goes first then Maker can choose $f$ such $\f(f)\leq (b+1)n$ and this is optimal. If Maker goes first then Maker can choose $f$ such $\f(f)\leq (n-1)b+n$ and this is optimal
\end{theorem}
We have the following corollary: Suppose that instead of paying $j,j=f(i)$, Maker pays $X(i,j)$ where the $X$-values are independent uniform $[0,1]$ random variables.
\begin{corollary}\label{cor1}
In this uniform random scenario, Maker can pay at most $\m_b(b+1+o(1))$, w.h.p., where $\m_b$ is the solution to $\tfrac{\log(b+1)-1}{b+1}=\m-1-\log\m$. (Note that $\m_b=1+O(\tfrac{\log^{1/2}b}{b^{1/2}})$.)
\end{corollary}
In preparation for more complicated scenarios, suppose that we replace $\N$ by $[m]$ i.e. we insist that $f(i)\leq m$.
\begin{theorem}\label{th2}
Suppose that $m>\sum_{i=1}^n1/i$. Then w.h.p. Maker can choose $f$ such $f(i)\leq m$ for $i\in [n]$ and $\f(f)\leq (b+1)n$ and this is optimal.
\end{theorem}
\begin{corollary}\label{cor2}
In this uniform random scenario, Maker can pay at most $\m_b(b+1+o(1))$, w.h.p.
\end{corollary}
Now we turn to a more complex problem. Here the {\em board} is the set of edges of the complete (loopless) digraph $\vK_n$. There is an $n\times n$ cost matrix $C$. For each $i\in [n]$ we have $C(i,i)=\infty$ and $C(i,j),j=1,2,\ldots [n]\setminus\set{i}$ is an independent  uniform random permutation of $[n-1]$. Maker's aim is to construct a spanning arborescence $T$ of low total cost, where the cost $C(T)=\sum_{(i,j)\in E(T)}C(i,j)$. (A spanning arborescence is a spanning tree whose edges have been oriented away from one vertex, the root.) 

After $t$ rounds, Maker will have selected a set of $t$ edges $M_t$ that induce a digraph where each vertex has out-degree at most one. Furthermore, these edges induce a forest when edge orientation is ignored. Similarly, Breaker will have selected a set $B_t$ of $bt$ edges where $M_t\cap B_t=\emptyset$.
\begin{theorem}\label{th3}
W.h.p., over the random choice of $C$,  Maker can construct an arborescence $T$ of cost $C(T)\leq \big(b / \th^*+b+1+o(1)\big)n$, where $\th^*\approx 0.2938...$. ($\th^*$ is the solution to $(1+\th^*)\log(1+\th^*)-\th^*\log\th^*=\log 2$.)
\end{theorem}
We doubt that this is optimal for Maker.

{\bf Conjecture:} W.h.p. Maker can construct an arborescence $T$ of cost $C(T)\leq (b+1+o(1))n$.

We base this conjecture on the fact that if $f:[n]\to[n]$ is a uniform random mapping then the digraph with edge-set $\set{(i,f(i))}$ is almost a random arborescence, i.e. it differs from a arborescence by $O(\log n)$ edges.
\begin{corollary}\label{cor3}
In this uniform random scenario, Maker can pay at most $\m_b\big(b / \th^*+b+1+o(1)\big)$, w.h.p.
\end{corollary}
We now consider the undirected versions of Theorem \ref{th3} and Corollary \ref{cor3}. I.e., now we have a weighted complete graph and Maker wishes to build a low cost spanning tree. We note that Hefetz, Kupferman, Lellouche and Vardi \cite{HKL} considered a worst-case version of this problem, in the context of finding a maximum weight spanning tree.

In the following theorem the edges of the complete graph are given independent uniform $[0,1]$ costs.
\begin{theorem}\label{th4}
Maker can build a spanning tree of cost at most $4\m_b\big(b / \th^*+b+1+o(1)\big)$, w.h.p.
\end{theorem}
If Maker and Breaker do not think too hard then they might come up with the following strategies for playing the Spanning Tree Game. Let us call the following strategies for Maker and Breaker the {\em greedy strategies}. The edges of $K_n$ are sorted into increasing order. Maker always chooses the cheapest edge available and Breaker always deletes the $b$ cheapest edges that Maker might want to use in the next round. 
\begin{theorem}\label{th5}
If Maker and Breaker employ greedy strategies then w.h.p. the cost of Maker's tree is asymptotically equal to $(b+1)\z(3)$.
\end{theorem}
\section{Mappings: proof of Theorem \ref{th1}}\label{sec1}
For $i\notin M_t$ we let $r_i(t)=\min\set{j:(i,j)\notin B_t}$ be the minimum of the possible values for $f(i)$ available to Maker at the end of round $t$. Let $\d_i(t)$ be the indicator for $i\notin M_t$. We use the following potential
\[
\F(t)=\sum_{i=1}^nr_i(t)\d_i(t).
\]
In round $t+1$, Maker's strategy is simply to choose an arbitrary $i\notin M_t$ and put $f(i)=r_i(t)$ and so that $M_{t+1}=M_t\cup \{i\}$. This reduces $\F$ by $r_i(t)$. On the other hand Breaker will add $b$ ordered pairs to $B_t$ and will be able to increase $\F$ by $b$. We get slightly different answers depending on who goes first. 

Suppose that Breaker goes first. Now $\F(0)=n$ and $\F(n)=0$. So we have
\[
-n=\sum_{i=0}^{n-1}(\F(t+1)-\F(t))=\sum_{i=0}^{n-1}(b-r_i(t))=nb-\sum_{i=0}^{n-1}r_i(t)=nb-\sum_{i=1}^nf(i).
\]
If Maker goes first then we replace $b-r_{n-1}(t)$ by $-r_{n-1}(t)$ because Breaker does not get an $n$th turn. This gives $\f(f)=(n-1)b+n$.

Both results are optimal, because any other choice for Maker with $f(i)>r_i(t)$ yields a higher value for $\f$. This completes the proof of Theorem \ref{th1}.

We complete this section with a proof of Corollary \ref{cor1}. We define permutations $\p_i$ such that $X(i,\p_i(j+1))\geq X(i,\p_i(j))$ for $1\leq j<n$. Then, if Maker were to choose $f(i)=j$ as in Theorem \ref{th1}, then Maker is charged $X(i,j)$. 

We need the following lemma from Frieze and Grimmett \cite{FG}.
\begin{lemma}\label{FG}
Suppose that $k_1+k_2+\cdots+k_M\leq aN$, and $Y_1,Y_2,\ldots,Y_M$ are independent random variables with $Y_i$ distributed as the $k_i$th minimum of $N$ independent uniform [0,1] random variables. If $\m>1$ then
\[
\Pr\left(Y_1+\cdots+Y_M\geq \frac{\m aN}{N+1}\right)\leq
e^{aN(1+\log \m-\m)}.
\]
(The lemma in \cite{FG} is given in terms of $a\log N$ instead of $aN$. There is nothing in the proof of that lemma that precludes us from replacing $a$ by $aN/\log N$.)
\end{lemma}
Now naively, we could observe that $\E(X(\p_i(j)))=j/(n+1)$ and then, at least in expectation, we could replace a cost $j$ in the function model, $j/(n+1)$, giving us a bound of $b+1$. Breaker however does affect the choice which $j$ Maker will choose and so we feel that we are forced to take a union bound over Maker's possibilities. This leads to the claimed inflated constant.

There are at most $\binom{(b+1)n-1}{n-1}$ choices for $\sum_{i=1}^nf(i)$ to add up to $(b+1)n$. Let $F$ denote this set of choices. For a fixed $\m>1$ we have
\mults{
\Pr\brac{\exists f\in F:\sum_{i=1}^nX(i,\p_i(f(i)))\geq \m(b+1)}\leq \binom{(b+1)n-1}{n-1}e^{-(b+1)n(1+\log\m-\m)}\\ \leq \Big((b+1)e^{-(1+(b+1)(1+\log\m-\m))}\Big)^n = o(1),
}
if $\m>\m_b$. This completes the proof Corollary \ref{cor1} and Corollaries \ref{cor2}, \ref{cor3} can be proved in the same manner. In particular note that the bound claimed in Corollary \ref{cor2} (resp. Corollary \ref{cor3}) is asymptotically just $\m_b$ times the coefficient of $n$ from Theorem \ref{th2} (resp. Theorem \ref{th3}).
\section{Mappings: proof of Theorem \ref{th2}}
Here Maker has to be more careful and exploit a version of the Box Game of Chv\'atal and Erd\H{o}s \cite{CE}. In this game there are $n$ disjoint sets $A_1,A_2,\ldots,A_n$, the boxes. In each round BoxMaker removes $p$ elements from $A=\bigcup_{i=1}^nA_i$ and BoxMaker removes $q$ boxes. In the context of our mapping game, we let $A_i=\set{(i,j),j\in[m]}$. Maker takes the role of BoxBreaker with $p=1$ and Breaker takes the role of BoxMaker with $q=b$. After $t$ rounds there will be $n-t$ boxes remaining and their contents will have been reduced. We will assume that BoxMaker (Breaker) goes first and the BoxBreaker (Maker) always chooses a remaining box $A_i$ of minimum size (and puts $f(i)$ equal to the smallest element left in $A_i$).

Theorem 3.4.1 of \cite{Book} shows that if $|A_i|=m>\sum_{i=1}^n1/i$ then BoxBreaker (Maker) has a strategy (described above) that garantees them a win. Thus given our lower bound on $m$, we see that Maker can finish the game. Also, by the analysis of Section \ref{sec1}, she will end with a value of $\f(f)\leq (b+1)n$. This completes the proof of Theorem \ref{th2}.
\section{Arborescence: proof of Theorem \ref{th3}}First we describe Maker's strategy. Let $F_t$ denote the set of oriented trees induced by $M_t$. Each of these trees/components will have a root and Maker must choose an edge leaving one of these roots.

On Maker's turn, if there is no root of any component on at most $n/2$ vertices such that Breaker has taken at least $n^\b$ edges from the root, we call this a {\em normal turn} for Maker. Here $0<\b<1$ is a constant to be determined later. On a normal turn, Maker always chooses the root $i$ of the smallest component $K$ (in case of a tie, say choose the least indexed root). Maker then chooses the edge $(i,j)\notin B_t$ that (i) minimises $C(i,k),k\notin B_t$ and (ii) does not point into $K$. We call this the {\em sensible} choice from root $i$.

If it is not a normal turn we say it is an {\em emergency turn}. On an emergency turn, Maker picks some component on at most $n/2$ vertices from which Breaker has removed at least $n^\b$ edges, and Maker makes the sensible choice out of this root. We will argue that w.h.p. this strategy gets Maker an arborescence consisting of a set of edges whose total cost is at most $(b/\th^*+b+1 +o(1))n$.

Consider the event to the contrary, i.e. that the total cost in the end is say $(b/\th^*+b+1 +3\e)n$ for some fixed $\e>0$, meaning that the extra cost paid for edges pointing within components is $(b/\th^*+3\e) n$, over and above the $(b+1)n$ achievable for Theorems \ref{th1}, \ref{th2}. We bound the probability of this event as follows. First it is at most $\p_1 + \p_2+ \p_3$, defined as follows. $\p_1$ is the probability that we pay an extra $\e n$ on normal steps falling in the first $n_0$ steps, where $n_0 = n-n^{\a}$ for some $0 < \a <1$. $\p_2$ is the probability of paying $\e n$ on normal steps after step $n_0$. $\p_3$ is the probability of paying $bn/\th^* + \e n$ on emergency steps. 

First we bound $\p_1$. Let $N_1$ be the set of normal steps in the first $n_0$ steps. Then we union bound over choices for numbers $a_t, t \in N_1$ where the intended meaning of the $a_t$ is as follows. At step $t$, Maker takes an edge from a root to another arborescense, paying an extra $a_t$ for edges pointing into its own component. So we have $\sum_t a_t = \e n$, and the number of choices for the $a_t$ is $\binom{\e n+|N_1|-1}{|N_1|-1}= \exp\{O(n)\}$. 

We will also union bound over sequences $x_t, t \in N_1$ which will mean the following. At step $i$, if $v$ is the root that Maker is taking an edge from and Maker chooses the edge of cost $r_t$, $x_t = r_t-1 - a_t$. In other words, the reason Maker has to choose the edge of cost $r_t = 1 + a_t + x_t$ is because among the edges of smaller cost, $a_t$ of them point into $v$'s component and $x_t$ edges have been taken by Breaker. In particular, knowing $x_t$ and $a_t$ tells us the cost $r_t$ of the edge that Maker will choose at step $t$. We have $\sum x_t \le bn$, and the number of choices for the $x_t$ is $\exp\{O(n)\}$. 

Having fixed the $a_t$ and $ x_t$ we will reveal the random digraph step by step as the game runs. More specifically, on Maker's turn at step $t\leq n-1$ we see the current component structure which uniquely determines the root, say $v$ which Maker will take an edge from. We reveal the costs of all edges coming from $v$, which determines which edge Maker will take (i.e. the lowest cost edge which is not taken by Breaker and which does not point into $v$'s component). Recall that among all the edges from $v$ not taken by Breaker, the $a_t$ lowest in cost all point into $v$'s component. The size of the smallest component is at most $\rdown{\frac{n}{n-t+1}}$ and so each out-edge has probability at most $\frac{\frac{n}{n-t+1}-1}{n}=\frac{t}{n(n-t+1)}$ of pointing into the component. There are at most $\binom{a_t + x_t}{a_t}$ choices for the costs of the $a_t$ edges pointing within the component. 
 Thus we bound the probability $\p_1$  by 
\begin{align}
\exp\{O(n)\} \prod_{t\in N_1} \binom{a_t + x_t}{a_t}\bfrac{t}{n(n-t+1)}^{a_t} &\le \exp\{O(n)\}   \prod_{t\in N_1} \bfrac{n-n^{\a}}{n^{1+\a}}^{a_t} \le n^{-\a \e n + o(n)}  = o(1). \label{eq1}
\end{align}

We bound $\p_2$ similarly. Let $N_2$ be the set of normal steps after step $n_0$. We choose numbers $a_t, t \in N_2$ adding up to $\e n$ and there are at most $\binom{\e n+n^\a}{n^\a} = \exp\{o(n)\}$ choices here. Likewise there are $\exp\{o(n)\}$ choices for the $x_t, i \in N_2.$ To bound the number of choices for the interior edges we can instead choose Breaker's edges and there are at most $n^{x_t}$ ways to do that.  So $\p_2$ is at most
\begin{align*}
\exp\{o(n)\} \prod_{t\in N_2}n^{x_t}\bfrac{t}{n(n-t+1)}^{a_t}&\leq \exp\{o(n)\} \times n^{n^{\a+\b}} \times 2^{-e n} = o(1)
\end{align*}
assuming only that $\a+\b<1$.

Finally we bound $\p_3$. Note first that there are at most $bn^{1-\b}$ emergency steps, so we union bound over at most $n^{bn^{1-\b}} = \exp\{o(n)\}$ choices for emergency steps. Let $N_3$ be the set of emergency steps. We fix numbers $a_t\geq n^\b, t \in N_3$ adding up to say $an$ where $a=b/\th^*+\e$, there being $\exp\{o(n)\}$ choices. We also fix numbers $x_t, t \in N_3$ adding up to say  $X \le b n$, there being $\exp\{o(n)\}$ choices. Since on an emergency step we always have a root of a component on at most $n/2$ vertices, the probability of an edge landing in the same component is at most $1/2$. Thus we bound $\p_3$ by
\beq{p3}{
\exp\{o(n)\} \prod_{t\in N_3}\binom{a_t + x_t}{a_t}\bfrac{1}{2}^{a_t} \le 2^{-an+o(n)} \prod_{t\in N_3}e^{(a_t+x_t)\log(a_t+x_t) - a_t \log(a_t) - x_t \log(x_t)}.
}
Suppose now that we fix the values for $a_t,t\in N_3$. Let 
\[
\f=\sum_{t\in N_3}f(a_t,x_t)\text{ where }f(a,x)=(a+x)\log(a+x) - a \log(a) - x \log(x).
\]
We argue next that to maximise $\f$ we must have $x_t/a_t$ taking the same value for all $t\in N_3$. Consider the function $g(x)=f(a,x)+f(b,L-x)$ for some $a,b,L>0$. Then 
\begin{align*}
g'(x)&=\log(a+x)-\log(x)-\log(b+L-x)+\log(L-x).\\
g''(x)&=-\frac{a}{x(a+x)}-\frac{b}{(L-x)(b+L-x)}<0.
\end{align*}
So, $g$ is strictly concave and its derivative vanishes when $(a+x)(L-x)=x(b+L-x)$ equivalently when $\tfrac{a}{x}=\tfrac{b}{L-x}$. Thus, if we fix $a_t,t\in N_3$ and maximise over $x_t,t\in N_3$ then we have $x_t=\th a_t$ for $t\in N_3$ where $\th=X/an\leq b/a$. We can therefore bound the product in \eqref{p3} by
\beq{p4}{
 \prod_{t\in N_3}e^{(a_t(1+\th))\log(a_t(1+\th)) - a_t \log(a_t) - \th a_t \log(\th a_t)}=\prod_{t\in N_3}e^{a_t((1+\th)\log(1+\th)-\th(\log\th))}=e^{an((1+\th)\log(1+\th)-\th\log\th)}.
}
Now recall that $\th^*\approx 0.2938...$ is defined as the root of $(1+\th^*)\log(1+\th^*)-\th^*\log\th^*=\log 2$. Since $a > b/\th^*$ we have $\th^* > b/a \ge \th$ and so $(1+\th)\log(1+\th)-\th\log\th < \log 2$. Now by \eqref{p3} and \eqref{p4} we have $\p_3 = o(1)$. 

We finally note that with the above Maker strategy, if a component reaches size greater than $n/2$ then its root will become the root of the final arborescence. This completes the proof of Theorem \ref{th3}.

\section{Spanning Trees 1: proof of Theorem \ref{th4}}
We reduce this Theorem \ref{th3}. We replace each edge $\{i,j\}$ with a pair of directed edges $(i,j),(j,i)$. Each directed edge $(i,j)$ is given a random cost $\wC(i,j)$, which is an independent copy of the $[0,1]$ random variable $Z$ where $\Pr(Z>x)=(1-x)^{1/2}\leq 1-\tfrac{x}2$ for $0\leq x\leq 1$. Then if $Z_1,Z_2$ are two independent copies of $Z$, then $\min\{Z_1,Z_2\}$ is distributed as a uniform $[0,1]$ random variable. This is a nice idea, employed by Walkup \cite{W} in bounding the expected value of a random assignment problem. Note that $Z$ is dominated by $2U[0,1]$.

Given the above construction, Maker builds a spanning arborescence. The factor 4 comes from two sources: first of all using $Z$ in place of $U[0,1]$ asymptotically doubles the cost of each selected edge. The value of the $k$th smallest order statistic of $n$ copies of $Z$ is at most doubled in expectation, due to the dominance claimed in the previous paragraph. The other source, is that if Breaker deletes $(i,j)$ then he must also delete $(j,i)$. So we double Breaker's power by replacing $b$ by $2b$. Of course, the factor $\m_b$ comes from using Lemma \ref{FG}, just as in the proof of Corollary \ref{cor1}.
\section{Spanning Trees 2: proof of Theorem \ref{th5}}
Suppose that the edges of $K_n$ are sorted into increasing order of weight, $e_1,e_2,\ldots,e_N,N=\binom{n}{2}$. Maker selects edges $e_{t_1},e_{t_2},\ldots,e_{t_{n-1}}$ and we approximate Maker's cost by $\frac1N\sum_{i=1}^{n-1}t_t$. Let $g_i=t_{i+1}-t_i$ and $T_k=t_1+t_2+\cdots +t_k$.

For $c>0$ we define $x=x(c)$ by (i) $0<x\leq 1$ and $xe^{-x}=ce^{-c}$. Let $m=cn/2$ then q.s.\footnote{A sequence of events $\cE_n,n\geq 1$ occurs {\em quite surely} q.s. if $\Pr(\neg\cE_n)=O(n^{-K})$ for an constant $K>0$.} $G_{n,m}$ contains a unique giant component of size $(1-\frac{x}c)n\pm n^{2/3}$. Let $\t(c)=\frac{1}c\sum_{k=1}^\infty\frac{k^{k-2}}{k!}(ce^{-c})^k$. Then q.s. $\k(G_{n,m})=t(c)n\pm n^{2/3}$. Then define
\[
t(k)=k,1\leq \frac{n}2\text{ and }t(k+1)=t(k)+\frac{1}{1-\brac{1-\frac{x_k}{c_k}}^2}\text{ for }k>\frac{n}2,
\]
where $\t(c_k)n=k$ and $0<x_k<1$ and $x_ke^{-x_k}=c_ke^{-c_k}$.

The analysis of the greedy algorithm in Frieze \cite{F} implies that w.h.p. $\sum_{k=1}^{n-1}t(k)\sim N\z(3)$.

Assume for now that Maker follows the greedy algorithm: Maker chooses $e_0$ in round zero and then in round $i,i\geq 1$ Maker chooses the cheapest edge that joins two distinct components. Breaker next chooses $b$ edges that would join two disitinct components in Maker's forest. 

Let $\cF_k$ denote the Maker's forest after the selection of $k$ edges, i.e. at the begining of round $k$. We note that $f_{k+1}$ is a randomly chosen edge that (i) avoids Breaker's edges and (ii) joins two distinct components of $\cF_k$. 

Assume that after $k\geq n/2$ rounds, $\cF_k$ has a typical structure i.e. a giant of size $\brac{1-\frac{x}c}n$ and all other components of size $O(\log n)$. Then the expected increase $t_{k+1}-t_k$ is asymptotically equal to $\frac{b+1}{1-\brac{1-\frac{x_k}{c_k}}^2}$ and the new forest component sizes remain ``typical''.

It is known (see for example Frieze and Karo\'nski \cite{FK}, Chapter 18) that if $L_n$ denotes the weight of the minimum spanning tree then 
\[
\E(L_n)=\int_{p=0}^1(\E(\k(G_n,p)-1)dp \sim \int_{c=0}^\infty \t(c)dc= \z(3).
\]
Now the effect of deleting $b$ edges is to slow the construction of the spanning tree and replace $t(c)$ by $t(c/(b+1))$ producing a tree of  expected weight asymptotically equal to 
\[
\int_{c=0}^\infty t(c/(b+1))dc=(b+1)\int_{c=0}^\infty \t(c)dc.
\]
\section{Final thoughts}
We have studied an interesting class of Maker-Breaker games where Maker's goal is build something cheaply. Our results are not all tight and we believe that there is a general meta theorem that states for such games, the existence of Breaker increases the cost of the optimum solution by a factor of $(b+1)$ on average.

\end{document}